\documentclass[11pt,fleqn]{article}
\usepackage{amssymb}
\usepackage{amsfonts}
\usepackage{amsmath}
\usepackage{amsthm}
\usepackage{graphicx}
\usepackage{epsfig}
\usepackage{psfrag}
\usepackage{color}
\usepackage{dsfont}
\makeatletter
\xdef\@endgadget#1{{\unskip\nobreak\hfil\penalty50\hskip1em\hbox{}\nobreak
    \hfil#1\parfillskip=0pt\finalhyphendemerits=0\par}}
\def\@qedsymbol{${}_\blacksquare$}
\def\qed{\@endgadget{\@qedsymbol}}
\newtheorem{lemma}{Lemma}[section]
\newtheorem{theorem}[lemma]{Theorem}

\newtheorem{example}[lemma]{Example}
\newtheorem{definition}[lemma]{Definition}

\newtheorem{proposition}[lemma]{Proposition}
\newtheorem{remark}[lemma]{Remark}
\newcommand{\mR}{\mathbb{R}}

\newcommand{\diag}{\mathrm{diag}}
\newcommand{\Exp}{\mathrm{Exp}}
\newcommand{\Ln}{\mathrm{Ln}}
\DeclareMathOperator{\im}{im} 
\DeclareMathOperator{\rank}{rank}

\DeclareMathOperator{\spa}{span}
\def\startmodif{\color{red}}
\def\stopmodif{\color{black}}

\def\BibTeX{{\rm B\kern-.05em{\sc i\kern-.025em b}\kern-.08em
    T\kern-.1667em\lower.7ex\hbox{E}\kern-.125emX}}

\title{\LARGE \bf On the Mathematical Structure of Balanced Chemical Reaction Networks Governed by Mass Action Kinetics
}

\author{Arjan van der Schaft\thanks{Johann Bernoulli Institute for Mathematics and Computer Science, University of Groningen, e-mail: {\tt\small A.J.van.der.Schaft@rug.nl}}\and Shodhan Rao\thanks{Systems Biology Centre for Energy Metabolism and Ageing, University of Groningen, email: {\tt \small s.rao@umcg.nl}}\and Bayu Jayawardhana\thanks{Discrete Technology and Production Automation, University of Groningen, email: {\tt \small b.jayawardhana@rug.nl}} 
}

\begin{document}

\maketitle
\thispagestyle{empty}
\pagestyle{empty}

\begin{abstract}
Motivated by recent progress on the interplay between graph theory, dynamics, and systems theory, we revisit the analysis of chemical reaction networks described by mass action kinetics. For reaction networks possessing a thermodynamic equilibrium we derive a compact formulation exhibiting at the same time the structure of the complex graph and the stoichiometry of the network, and which admits a direct thermodynamical interpretation. This formulation allows us to easily characterize the set of equilibria and their stability properties. Furthermore, we develop a framework for interconnection of chemical reaction networks. Finally we discuss how the established framework leads to a new approach for model reduction.
\end{abstract}

\maketitle              

\section{Introduction}
Large-scale chemical reaction networks arise abundantly in bio-engineering and systems biology. A very simple example involving only two chemical reactions in the glycolytic pathway is given by the following two coupled reactions
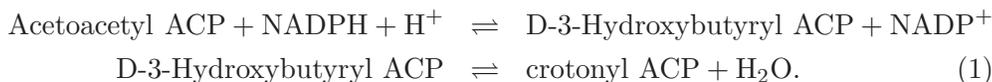
\begin{eqnarray}\label{eq:eg}
\nonumber \text{Acetoacetyl ACP}+ \text{NADPH} +\text{H}^+  & \rightleftharpoons & \text{D-3-Hydroxybutyryl ACP}+\text{NADP}^+\\
\text{D-3-Hydroxybutyryl ACP} & \rightleftharpoons & \text{crotonyl ACP} +\text{H}_2\text{O}.
\end{eqnarray}
The compounds Acetoacetyl ACP, NADPH, $\text{H}^+$, D-3-Hydroxybutyryl ACP, NADP, crotonyl ACP and $\text{H}_2\text{O}$ involved in the reaction network are called the {\it chemical species} of the network. The most basic law prescribing the dynamics of the concentrations of the various species involved in a chemical reaction network is the \emph{law of mass action kinetics}, leading to polynomial differential equations for the dynamics of each species. Large-scale networks thus lead to a high-dimensional set of coupled polynomial differential equations, which are usually difficult to analyze.

In order to gain insight into the dynamical properties of complex reaction networks it is important to identify their underlying mathematical structure, and to express their dynamics in the most compact way. In line with the recent surge of interest in network dynamics at least two aspects should be fundamental in such a mathematical formulation: (1) a graph representation, and (2) a specific form of the differential equations.

The graph representation of chemical reaction networks is not immediate, since a chemical reaction (the obvious candidate for identification with the edges of a graph) generally involves more than two chemical species (the most basic candidate for identification with the vertices). We will follow an approach that has been initiated and developed in the work of Horn \& Jackson \cite{HornJackson, Horn} and Feinberg, see e.g. \cite{Feinberg, Feinberg1}, by associating the {\it complexes} of the chemical reaction network (i.e., the left- and right-hand sides of the reactions) with the vertices of the graph\footnote{For alternative approaches based on species-reaction graphs and Petri-nets we refer to e.g. \cite{craciun, Angeli, AngeliEJC}.}. The resulting directed graph, called the {\it complex graph} in this paper, is characterized by its incidence matrix. The expression of the complexes in the chemical species defines an extended stoichiometric matrix, called the {\it complex stoichiometric matrix}, which is immediately related to the standard stoichiometric matrix through the incidence matrix of the complex graph; see e.g. \cite{Othmer}, \cite{AngeliEJC}.

In order to derive a specific form of differential equations we will start with the basic assumption that the chemical reaction rates are governed by {\it mass action kinetics}. As an initial step we then derive a compact form of the dynamics involving a non-symmetric weighted Laplacian matrix of the complex graph. The basic form of these equations can be already found in the innovative paper by Sontag \cite{Sontag}. The main part of the paper is however devoted to a subclass of mass action kinetics chemical reaction networks, where we assume the existence of a thermodynamical equilibrium, or equivalently, where the {\it detailed balance equations} are assumed to admit a solution. We will call such chemical reaction networks {\it balanced chemical reaction networks}. Balanced chemical reaction networks are necessarily reversible but involve additional conditions on the forward and reverse reaction rate constants (usually referred to as the Wegscheider conditions; see \cite{FeinbergWegscheider}).
For such balanced chemical reaction networks we will be able to derive a particularly elegant form of the dynamics, involving a {\it symmetric weighted Laplacian matrix} of the complex graph. Furthermore, it turns out that this form has a direct thermodynamical interpretation, and in fact can be regarded as a graph-theoretic version of the formulation derived in the work of Katchalsky, Oster \& Perelson \cite{op, opk}.

The obtained form of the equations of balanced chemical reaction networks will be used to give, in a very simple and insightful way, a characterization of the set of equilibria, and a proof of the asymptotic convergence to a unique thermodynamic equilibrium corresponding to the initial condition of the system. Similar results for a different class of mass action kinetics reaction networks, in particular {\it weakly reversible} networks with zero-deficiency, have been derived in the fundamental work of Horn \cite{Horn} and Feinberg \cite{Feinberg, Feinberg1}, which was an indispensable source of concepts and tools for the work reported in this paper. 

Subsequently we show how this form of the dynamics of chemical reaction networks can be extended to {\it open} reaction networks; i.e., networks involving an influx or efflux of some of the chemical species, called the {\it boundary chemical species}. Furthermore, we show how corresponding outputs can be defined as the chemical potentials of these boundary chemical species and how this leads to a physical theory of {\it interconnection} of open reaction networks, continuing upon the work of Oster \& Perelson \cite{opk, opcircuits}.

In the final part of the paper we make some initial steps in showing how the derived form of balanced reaction networks can be utilized to derive in a systematic way {\it reduced-order models}. This reduction procedure is based on the Laplacian matrix of the complex graph, and draws inspiration from a similar technique in electrical circuit theory, sometimes referred to as {\it Kron reduction}. The model reduction procedure leads to an approximating reduced-order model which is again a mass action kinetics balanced chemical reaction network.

Interestingly enough, the derived form of the equations of chemical reaction networks involving a weighted Laplacian matrix of the complex graph is reminiscent of the dynamics of the standard {\it consensus algorithm of multi-agent systems}, see e.g. \cite{Murray}. In fact, the form of the stability analysis in both cases is quite similar. Furthermore, 
the relation between the initial concentration of chemical species and the unique final equilibrium concentration is closely related to the recent concept of $\chi$-consensus in \cite{Cortes}.

\smallskip

The paper is organized as follows. In Section \ref{sec:CRNS}, we summarize the mathematical structure of a network of chemical reactions described in e.g. \cite{HornJackson, Horn, Feinberg, Feinberg1, Othmer}; see also \cite{AngeliEJC} for a clear account. In Section \ref{sec:MAK1}, we recall the law of mass action kinetics, and derive a general graph-theoretic formulation which is basically already contained in the innovative paper by Sontag \cite{Sontag}. Next, for the important subclass of chemical reaction networks admitting a thermodynamic equilibrium ({\it balanced  reaction networks}) we derive a particularly simple formulation involving a symmetric weighted Laplacian matrix for the complex graph. Before providing the thermodynamical interpretation of this formulation, we discuss the graph-theoretic treatment of the complex graph, largely following the exposition in \cite{Othmer}, and its implications for the established formulation of balanced reaction networks. In Section \ref{sec:ICSM} we utilize the developed formulation to derive a simple characterization of the set of equilibria, and we give a Lyapunov analysis to show the convergence to a unique equilibrium depending on the initial concentration vector. At the end of this section we comment on the relation with consensus algorithms. In Section \ref{sec:BF} the framework is extended to reaction networks with inflows and outflows, and we discuss interconnection of such reaction networks through shared boundary chemical species. The initial steps in the application of the established formulation of balanced reaction networks to model reduction are discussed in Section \ref{sec:modelreduction}. Finally, Section \ref{sec:Conc} contains the conclusions, and some topics of ongoing research.

\bigskip

\noindent\emph{\bf Notation}:  The space of ${n}$ dimensional real vectors is denoted by $\mathbb{R}^{{n}}$,
and the space of ${m}\times {n}$ real matrices by $\mathbb{R}^{{m}\times {n}}$. The space of ${n}$ dimensional real vectors consisting of all strictly positive entries is denoted by $\mR_+^{n}$ and the space of ${n}$ dimensional real vectors consisting of all nonnegative entries is denoted by $\bar{\mR}_+^{n}$. The rank of a real matrix $A$ is denoted by $\rank A$. Given $a_1,\ldots,a_n \in \mR$, $\mbox{diag}(a_1,\ldots,a_n)$ denotes the diagonal matrix with diagonal entries $a_1,\ldots,a_n$; this notation is extended to the block-diagonal case when $a_1,\ldots,a_n$ are real square matrices. Furthermore, $\ker A$ and $\spa A$ denote the kernel and span respectively of a real matrix $A$. If $U$ denotes a linear subspace of $\mathbb{R}^m$, then $U^{\perp}$ denotes its orthogonal subspace (with respect to the standard Euclidian inner product). $\mathds{1}_m$ denotes a vector of dimension $m$ with all entries equal to 1. The time-derivative $\frac{dx}{dt}(t)$ of a vector $x$ depending on time $t$ will be usually denoted by $\dot{x}$. 

Define the mapping
$\mathrm{Ln} : \mathbb{R}_+^m \to \mathbb{R}^m, \quad x \mapsto \mathrm{Ln}(x),$
as the mapping whose $i$-th component is given as
$\left(\mathrm{Ln}(x)\right)_i := \mathrm{ln}(x_i).$
Similarly, define the mapping
$\mathrm{Exp} : \mathbb{R}^m \to \mathbb{R}_+^m, \quad x \mapsto \mathrm{Exp}(x),$
as the mapping whose $i$-th component is given as
$\left(\mathrm{Exp}(x)\right)_i := \mathrm{exp}(x_i).$
Also, define for any vectors $x,z \in \mathbb{R}^m$ the vector $x \cdot z \in \mathbb{R}^m$ as the element-wise product $\left(x \cdot z\right)_i :=x_iz_i, \, i=1,2,\ldots,m,$ and the vector $\frac{x}{z} \in \mathbb{R}^m$ as the element-wise quotient $\left(\frac{x}{z}\right)_i := \frac{x_i}{z_i}, \, i=1,\cdots,m$. Note that with these notations $\Exp (x + z) = \Exp (x) \cdot \Exp (z)$ and $\Ln (x \cdot z) = \Ln (x) + \Ln (z), \Ln \left(\frac{x}{z}\right) = \Ln (x) - \Ln (z)$.

\section{Chemical reaction network structure} \label{sec:CRNS}
In this section we will survey the basic topological structure of chemical reaction networks. First step is the stoichiometry expressing the conservation laws of chemical reactions. A next innovative step was taken in the work of Horn \& Jackson and Feinberg \cite{HornJackson, Horn, Feinberg, Feinberg1} by defining the complexes of a reaction to be the vertices of a graph. We will summarize these achievements in a slightly more abstract manner, also making use of the exposition given in \cite{Othmer}; see also \cite{AngeliEJC} for a nice account.

\subsection{Stoichiometry}
Consider a chemical reaction network involving $m$ chemical species (metabolites), among which $r$ chemical reactions take place. The basic structure underlying the dynamics of the vector $x \in \bar{\mathbb{R}}_+^m$ of concentrations $x_i, i=1,\cdots, m,$ of the chemical species is given by the {\it balance laws}
\begin{equation}\label{balance}
\dot{x} = Sv,
\end{equation}
where $S$ is an $m \times r$ matrix, called the {\it stoichiometric matrix}. The elements of the vector $v \in \mathbb{R}^r$ are commonly called the (reaction) {\it fluxes}. The stoichiometric matrix $S$, which consists of (positive and negative) integer elements, captures the basic conservation laws of the reactions. For example, the stoichiometric matrix of the two coupled reactions involving the chemical species $X_1,X_2,X_3$ given as
\[
X_1 + 2X_2 \rightleftharpoons X_3 \rightleftharpoons 2X_1 + X_2
\]
is
\[
S = \begin{bmatrix} -1 & 2 \\ -2 & 1 \\ 1 & -1 \end{bmatrix},
\]
while the first reversible reaction in (\ref{eq:eg})
\[
\text{Acetoacetyl ACP}+ \text{NADPH} +\text{H}^+  \rightleftharpoons  \text{D-3-Hydroxybutyryl ACP}+\text{NADP}^+
\]
involving the species \{Acetoacetyl ACP, NADPH, $\text{H}^+$, D-3-Hydroxybutyryl ACP, $\text{NADP}^+$, crotonyl ACP, $\text{H}_2$O\} has the stoichiometric matrix
\[
S= \begin{bmatrix} -1 & -1 & -1 & 1 & 1 & 0 \end{bmatrix}^T
\]
In many cases of interest, especially in biochemical reaction networks, chemical reaction networks are intrinsically {\it open}, in the sense that there is a continuous exchange with the environment (in particular, flow of chemical species and connection to other reaction networks). This will be modeled by a vector $v_b \in \mathbb{R}^b$ consisting of $b$ {\it boundary} {\it fluxes} $v_b$, leading to an extended model
\begin{equation}\label{stoichiometry}
\dot{x} = Sv+ S_b v_b
\end{equation}
Here the matrix $S_b$ consists of mutually different columns whose elements are all $0$ except for one element equal to $1$ or $-1$. The non-zero elements correspond to boundary fluxes for part of the chemical species, called the {\it boundary chemical species}, with $+1$ denoting an uptake flux and $-1$ a demand flux. 

\begin{example}
The open reaction network
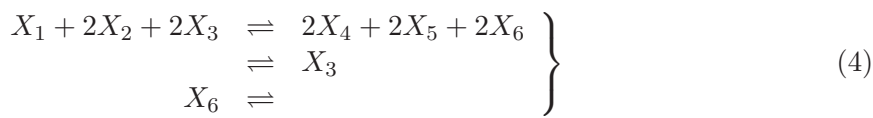
\begin{equation}\label{boundaryfluxexample}
\left.
\begin{array}{rcl}
X_1+2X_2 + 2X_3 & \rightleftharpoons & 2X_4+2X_5+2X_6 \\
 & \rightleftharpoons & X_3 \\
X_6 & \rightleftharpoons &
\end{array}\right\}
\end{equation}
has the extended stoichiometric matrix
\[
\begin{bmatrix} S & S_b \end{bmatrix} = \begin{bmatrix} -1 & | & 0 & 0 \\ -2 & | & 0 & 0 \\ -2 & | & 1 & 0 \\
2 & | & 0 & 0 \\ 2 & | & 0 & 0 \\ 2 & | & 0 & -1 
\end{bmatrix}
\]
with boundary chemical species $X_3$ and $X_6$.
\end{example}

Although the stoichiometry does not fully determine the dynamics of the chemical reaction network (for this to be the case the vector of fluxes $v$ needs to be expressed as a {\it reaction rate} $v=v(x)$ ) it already contains useful information about the network dynamics, {\it independent} of the precise form of the reaction rate $v(x)$. For example, if there exists an $m$-dimensional row-vector $k$ such that
\[
kS = 0
\]
then the quantity $kx$ is a {\it conserved quantity} for the dynamics $\dot{x} = Sv(x)$ for {\it all} possible reaction rates $v=v(x)$. Indeed, $\frac{d}{dt} kx = kSv(x) =0$. If $k \in \bar{\mathbb{R}}_+^m$, then the quantity $kx$ is commonly called a {\it conserved moiety}. Geometrically, for all possible fluxes the solutions of the differential equations $\dot{x} = Sv(x)$ starting from an initial state $x_0$ will always remain within the affine space $x_0 + \im S$.
In case of an open chemical reaction network $\frac{d}{dt} kx =0$ modifies into
\[
\frac{d}{dt} kx = kS_bv_b,
\]
expressing that the time-evolution of the quantity $kx$ only depends on the boundary fluxes $v_b$.

\subsection{The complex graph}
The network structure of a chemical reaction network cannot be directly captured by a graph involving the chemical species (since generally there are more than two species involved in a reaction). Instead, we will follow an approach originating in the work of Horn \& Jackson \cite{HornJackson, Horn} and Feinberg \cite{Feinberg, Feinberg1}, introducing the space of {\it complexes}. The set of complexes of a chemical reaction network is simply defined as the union of all the different left- and right-hand sides (substrates and products) of the reactions in the network. 

For example, the reaction network in (\ref{eq:eg}) entails four complexes, namely the substrates  Acetoacetyl ACP+ NADPH +$\text{H}^+$ and D-3-Hydroxybutyryl ACP+$\text{NADP}^+$, and the products D-3-Hydroxybutyryl ACP  and crotonyl ACP +$\text{H}_2$O. In general, a product complex of one reaction may be the substrate complex of another. Thus, the complexes corresponding to the two reactions
\[
2X_1 \rightleftharpoons X_1 + 2X_2 \rightleftharpoons X_2 + X_3
\]
are $2X_1, X_1 + X_2$ and $X_2 + X_3$. A complex may also be the product/substrate complex of more than one reaction.

The expression of the complexes in terms of the chemical species is formalized by an $m \times c$ matrix $Z$, whose $\rho$-th column captures the expression of the $\rho$-th complex in the $m$ chemical species. For example, the column expressing the complex $X_1+2X_2+2X_3$ in the chemical reactions (\ref{boundaryfluxexample}) is 
\[
\begin{bmatrix} 1 & 2 & 2 & 0 & 0 & 0 \end{bmatrix}^T 
\]
Note that by definition all elements of the matrix $Z$ are non-negative integers. 

Since the complexes are left- and right-hand sides of the reactions they can be naturally associated with the vertices of a {\it directed graph}, with edges corresponding to the reactions. The complexes on the left-hand side of the reactions are called the {\it substrate} complexes and those on the right-hand side of the reactions are called the {\it product} complexes. Formally, the reaction $ \sigma  \rightleftharpoons \pi$ between the $\sigma$-th and the $\pi$-th complex defines a directed edge with tail vertex being the $\sigma$-th complex and head vertex being the $\pi$-th complex. The resulting graph is called the {\it complex graph}. 

Recall, see e.g. \cite{Bollobas}, that any directed graph is defined by its {\it incidence matrix} $B$. This is an $c \times r$ matrix, $c$ being the
number of vertices and $r$ being the number of edges, with $(\rho,j)$-th
element equal to $-1$ if vertex $\rho$ is the tail vertex of edge $j$ and $1$ if vertex $\rho$ is the head vertex of edge $j$, while $0$ otherwise. 
Obviously there is a close relation between the matrix $Z$ and the stoichiometric matrix $S$, which is expressed as
\begin{equation}
S = ZB,
\end{equation} 
with $B$ the incidence matrix of the complex graph. For this reason we will call $Z$ the {\it complex stoichiometric matrix}. Hence the relation $\dot{x} = Sv$ between the fluxes $v$ through the chemical reaction network and the time-derivative of vector of chemical species $x$, cf. (\ref{balance}), can be also written as
\begin{equation}
\dot{x} = ZBv,
\end{equation}
with the vector $Bv$ belonging to the space of complexes $\mathbb{R}^c$.


\section{The dynamics of mass action kinetics chemical reaction networks}\label{sec:MAK1}
In this section we will derive a compact form for the dynamics of a chemical reaction network, whose reactions are described by {\it mass action kinetics}. After deriving a general form in Section 3.1 we will focus in the rest of the section on {\it balanced} chemical reaction networks. The resulting form of the dynamical equations for balanced reaction networks will be fundamental to the analysis of this dynamics in the subsequent sections.

\subsection{The general form of mass action kinetics}
The dynamics of the concentration vector $x$ (or equivalently, in case of a fixed volume, the vector $n$ of mole numbers) is given once the internal fluxes $v$ are specified as a function $v=v(x)$ of $x$, defining the {\it reaction rates}. The most basic model for specifying the reaction rates is {\it mass action kinetics}, defined as follows. Consider the single reaction
\[
X_1 + X_2 \rightleftharpoons X_3,
\]
involving the three chemical species $X_1,X_2,X_3$ with concentrations $x_1,x_2,x_3$. In mass action kinetics the reaction is considered to be a combination of the {\it forward reaction} $X_1 + X_2 \rightharpoonup X_3$ with forward rate equation $v^{\text{forw}}(x_1,x_2) = k^{\text{forw}}x_1x_2$ and the {\it reverse reaction}
$
X_1 + X_2 \leftharpoondown X_3,
$
with rate equation $v^{\text{rev}}(x_3) = k^{\text{rev}}x_3$. The constants $k^{\text{forw}}$, $k^{\text{rev}}$ are called respectively the {\it forward} and {\it reverse reaction constants}.
The net reaction rate given by mass action kinetics is thus
\[
v(x_1,x_2,x_3) = v^{\text{forw}}(x_1,x_2) - v^{\text{rev}}(x_3) = k^{\text{forw}} x_1x_2 - k^{\text{rev}}x_3
\]
In general, the mass action reaction rate of the $j$-th reaction of a chemical reaction network, from a substrate complex $\mathcal{S}_j$ to a product complex $\mathcal{P}_j$, is given as
\begin{equation}\label{eq:massaction}
v_j(x) = k_j^{\text{forw}} \prod_{i=1}^m x_i ^{Z_{i \mathcal{S}_j}} - k_j^{\text{rev}} \prod_{i=1}^m x_i^{Z_{i \mathcal{P}_j}},
\end{equation}
where $Z_{i \rho}$ is the $(i,\rho)$-th element of the complex stoichiometric matrix $Z$, and $k_j^{\text{forw}}, k_j^{\text{rev}} \geq 0$ are the forward/reverse reaction constants of the $j$-th reaction, respectively. Without loss of generality we will throughout assume that for every $j$ the constants $k_j^{\text{forw}}, k_j^{\text{rev}}$ are {\it not both equal to zero} (since in this case the $j$-th reaction is not active). 

This can be rewritten in the following way. Let $Z_{\mathcal{S}_j}$ and $Z_{\mathcal{P}_j}$ denote the columns of the complex stoichiometry matrix $Z$ corresponding to the substrate and the product complexes of the $j$-th reaction. Using the mapping $\mathrm{Ln}  : \mathbb{R}^c_+ \to \mathbb{R}^c$ as defined at the end of the Introduction, the mass action reaction equation (\ref{eq:massaction}) for the $j$-th reaction from substrate complex $\mathcal{S}_j$ to product complex $\mathcal{P}_j$ can be rewritten as
\begin{equation}\label{eq:massactionwithZ}
v_j(x)=k^{\text{forw}}_{j} \exp\big(Z_{\mathcal{S}_j}^T \mathrm{Ln}(x)\big)-k^{\text{rev}}_{j} \exp\big(Z_{\mathcal{P}_j}^T \mathrm{Ln}(x)\big).
\end{equation}
Based on the formulation in (\ref{eq:massactionwithZ}), we can describe the complete reaction network dynamics as follows. Let the mass action rate for the complete set of reactions be given by the vector $v(x)= \begin{bmatrix} v_1(x) & \cdots & v_r(x) \end{bmatrix}^T$. 
For every $\sigma,\pi \in \{1,\cdots, c\}$, denote by $a_{\sigma \pi}=k^{\text{rev}}_{j}$, $a_{\pi \sigma}=k^{\text{forw}}_{j}$ if $(\sigma,\pi)=(\mathcal{S}_j,\mathcal{P}_j)$, $j \in \{1,\ldots,r\}$ and $a_{\sigma \pi}=0$ elsewhere. Define the {\it weighted adjacency matrix} $A$ of the complex graph\footnote{Strictly speaking, $A$ is the adjacency matrix {\it not} of the complex graph, but of a directed {\it augmented} graph. Each $j$-th edge of the complex graph is replaced by {\it two} directed edges: one corresponding to the forward reaction with weight $k_j^{\text{forw}}$ and one corresponding to the reverse reaction with weight $k_j^{\text{rev}}$.} as the matrix with $(\sigma,\pi)$-th element $a_{\sigma \pi}$, where $\sigma,\pi \in \{1, \cdots,c\}$.
Furthermore, define the {\it weighted Laplacian} matrix  $L$ as the $c \times c$ matrix 
\begin{equation}\label{Laplacian}
L := \Delta - A
\end{equation}
where $\Delta$ is the diagonal matrix whose $(\rho,\rho)$-th element is equal to the sum of the elements of the $\rho$-th column of $A$.
Then it can be verified that the vector $Bv(x)$ for the mass action reaction rate vector $v(x)$ is equal to
\[
Bv(x) = -L \Exp \left(Z^T \Ln(x) \right),
\]
where the mapping $\Exp : \mathbb{R}^c \to \mathbb{R}^c_+$ has been defined at the end of the Introduction.
Hence the dynamics can be compactly written as
\begin{equation}\label{sontag}
\dot{x} = - Z L \mathrm{Exp} \left(Z^T \mathrm{Ln}(x)\right)
\end{equation}
A similar expression of the dynamics corresponding to mass action kinetics, in less explicit form, was already obtained in \cite{Sontag}.

\subsection{Balanced mass action kinetics}
In the rest of this paper we will focus on an important subclass of mass action chemical reaction networks. We need the following basic notions.
\begin{definition}
A vector of concentrations $x^* \in \mathbb{R}^m_+$ is called an {\it equilibrium}\footnote{Sometimes called a {\it kinetic} equilibrium, in order to distinguish it from a thermodynamical equilibrium.} for the dynamics $\dot{x} =Sv(x)$ if 
\[
Sv(x^*) = 0
\]
Furthermore, $x^* \in \mathbb{R}^m_+$ is called a {\it thermodynamic equilibrium} if 
\[
v(x^*)=0
\]
\end{definition}
Clearly, any thermodynamic equilibrium is an equilibrium, but not necessarily the other way around (since in general $S=ZB$ is not injective). 

\subsubsection{The existence of thermodynamic equilibria}
Necessary and sufficient conditions for the existence of a thermodynamic equilibrium can be derived in the following linear-algebraic way following \cite{FeinbergWegscheider}. These conditions are usually referred to as the {\it Wegscheider conditions}, generalizing the classical results of \cite{wegscheider}. 

Consider the $j$-th reaction from substrate $\mathcal{S}_{j}$ to product $\mathcal{P}_{j}$, described by the mass action rate equation \big(see (\ref{eq:massactionwithZ})\big)
\[
v_{j}(x) = k^{\text{forw}}_{j} \exp\big(Z_{\mathcal{S}_{j}}^T \mathrm{Ln}(x)\big)-k^{\text{rev}}_{j} \exp\big(Z_{\mathcal{P}_{j}}^T \mathrm{Ln}(x)\big)
\]
Then $x^* \in \mathbb{R}^m_+$ is a \startmodif thermodynamic \stopmodif equilibrium, i.e., $v(x^*)=0$, if and only if \startmodif
\begin{equation}\label{detailedbalance}
k^{\text{forw}}_{j} \exp\big(Z_{\mathcal{S}_{j}}^T \mathrm{Ln}(x^*)\big) = k^{\text{rev}}_{j} \exp\big(Z_{\mathcal{P}_{j}}^T \mathrm{Ln}(x^*)\big), \quad j=1, \cdots,r
\end{equation} \stopmodif
The equations (\ref{detailedbalance}), sometimes referred to as the {\it detailed balance equations}, can be rewritten as follows.
Define the {\it equilibrium constant} $K_{j}^{eq}$ of the $j$-th reaction as (assuming $k^{\text{rev}}_{j} \neq 0$)
\begin{equation}\label{eqconst}
K_{j}^{eq}: = \frac{k^{\text{forw}}_{j}}{k^{\text{rev}}_{j}}
\end{equation}
Then the detailed balance equations (\ref{detailedbalance}) are seen to be equivalent to
\begin{equation}\label{detailed}
K_{j}^{eq} = \exp \left( Z_{\mathcal{P}_{j}}^T \Ln(x^*) - Z_{\mathcal{S}_{j}}^T \Ln(x^*) \right), \quad j =1, \cdots,r
\end{equation}
Collecting all reactions, and making use of the incidence matrix $B$ of the complex graph, this amounts to the vector equation
\begin{equation}\label{wegscheider}
K^{eq} = \Exp \left(B^TZ^T \Ln (x^*) \right) = \Exp \left(S^T \Ln (x^*) \right),
\end{equation}
where $K^{eq}$ is the $r$-dimensional vector with $j$-th element $K_{j}^{eq}, j =1, \cdots,r$.

From here it is easy to characterize the existence of a thermodynamic equilibrium.
\begin{proposition}\label{thermoequilibrium}
There exists a thermodynamic equilibrium $x^* \in \mathbb{R}^m_+$ if and only if $ k_{j}^{\text{forw}} > 0, k_{j}^{\text{rev}} > 0,$ for all $j=1, \cdots,r$, and furthermore
\begin{equation}\label{wegscheider1}
\Ln \left(K^{eq}\right) \in \im S^T
\end{equation}
\end{proposition}
\begin{proof}
Clearly $ k_{j}^{\text{forw}} > 0, k_{j}^{\text{rev}} > 0, j=1, \cdots,r,$ is a necessary condition for the existence of a thermodynamic equilibrium. The existence of a vector $\Ln (x^*), x^* \in \mathbb{R}^m_+,$ satisfying (\ref{wegscheider}) is obviously equivalent to (\ref{wegscheider1}).
\end{proof}
It also follows that once a thermodynamic equilibrium $x^*$ is given, the set of {\it all} thermodynamic equilibria is described as follows.
\begin{proposition}\label{thermoequilibrium1}
Let $x^* \in \mathbb{R}^m_+$ be a thermodynamic equilibrium, then the set of {\it all} thermodynamic equilibria is given by 
\begin{equation}\label{equilibria}
\mathcal{E} := \{ x^{**} \in \mathbb{R}^m_+ \mid S^T \Ln \left(x^{**}\right) = S^T \Ln \left(x^{*}\right) \}
\end{equation}
\end{proposition}
For ease of exposition we will henceforth refer to mass action chemical networks possessing a thermodynamic equilibrium $x^* \in \mathbb{R}^m_+$ as {\it balanced mass action chemical networks}:
\begin{definition}
A chemical reaction network $\dot{x}=Sv(x)$ governed by mass action kinetics is called {\it balanced} if there exists a thermodynamic equilibrium $x^* \in \mathbb{R}^m_+$.
\end{definition}
As stated in Proposition \ref{thermoequilibrium}, a necessary condition for the existence of a thermodynamic equilibrium is the fact that for all reactions the forward and reverse reaction constants are both strictly positive. Thus all reactions of a balanced reaction network need to be at least {\it reversible}. 

\begin{remark}
Usually the conditions (\ref{wegscheider1}) are rewritten in the following more constructive form. By basic linear algebra (\ref{wegscheider1}) is satisfied if and only if for all row-vectors $\sigma$ satisfying $\sigma S^T = 0$ we have
\begin{equation}
\sum_{j =1}^r \sigma_{j} \ln K_{j}^{eq} = 0
\end{equation}
Putting back in the definition of the equilibrium constants (\ref{eqconst}) this is seen to be equivalent to the usual form of the Wegscheider conditions\footnote{Note that (\ref{wegscheider2}) only needs to be checked for a maximally independent set of row-vectors $\sigma$ satisfying $\sigma S^T = 0$. Furthermore, by writing this as $\sigma S^T = \sigma B^TZ^T =0$ this can be related to the topological structure of the complex graph and to the deficiency of the network; cf. \cite{FeinbergWegscheider}. Furthermore, the Wegscheider conditions admit a thermodynamical interpretation, since they are intimately related to the {\it passivity} properties of the network and the independence of the increase of Gibbs' energy of the path in the reaction network; see also \cite{opk,op}.}
\begin{equation}\label{wegscheider2}
\prod_{j =1}^r \left(k^{\text{forw}}_{j} \right)^{\sigma_{j}} = \prod_{j =1}^r \left(k^{\text{rev}}_{j}\right)^{\sigma_{j}}.
\end{equation}
\end{remark}

\subsubsection{The standard form of balanced mass action reaction networks}
For balanced mass action chemical reaction networks we can further rewrite the dynamics (\ref{sontag}) in the following useful way. Let $x^* \in \mathbb{R}_+^m$ be a thermodynamic equilibrium, i.e., $v(x^*)=0$. Consider the rewritten form (\ref{detailed}) of the `detailed balance' equations. These equations allow us to define the {\it balanced reaction constant}\footnote{Although these constants appear in the literature, see e.g. \cite{opk,op}, we could not find a standard name for them. Note that $\kappa_{j}(x^*) > 0 $ is defined as the common value of the forward and reverse reaction rate of the $j$-th reaction at thermodynamic equilibrium.} of the $j$-th reaction as
\begin{equation}\label{kappa}
\kappa_{j}(x^*) :=k_{j}^{\text{forw}} \exp \left(Z_{\mathcal{S}_{j}}^T \Ln (x^*)\right)=k_{j}^{\text{rev}} \exp \left(Z_{\mathcal{P}_{j}}^T \Ln (x^*)\right), \quad j = 1, \cdots, r
\end{equation}
Then the mass action reaction rate of the $j$-th reaction can be written as
\[
v_{j}(x)=\kappa_{j}(x^*)\left[\exp\left(Z_{\mathcal{S}_{j}}^T \mathrm{Ln}\left(\frac{x}{x^*}\right)\right)-\exp\left(Z_{\mathcal{P}_{j}}^T \mathrm{Ln}\left(\frac{x}{x^*}\right)\right)\right],
\]
where for any vectors $x,z \in \mathbb{R}^m$ the quotient vector $\frac{x}{z} \in \mathbb{R}^m$ is defined elementwise (see the end of the Introduction).
Defining the $r \times r$ diagonal matrix of balanced reaction constants as
\begin{equation}\label{K}
\mathcal{K}(x^*) := \mathrm{diag} \big( \kappa_1(x^*), \cdots, \kappa_r(x^*) \big)
\end{equation}
it follows that the mass action reaction rate vector of a balanced reaction network can be written as
\[
v(x) = - \mathcal{K}(x^*) B^T \mathrm{Exp} \left(Z^T \mathrm{Ln}\left(\frac{x}{x^*}\right)\right),
\]
and thus the dynamics of a balanced reaction network takes the form
\begin{equation}\label{masterequation}
\dot{x} = - Z B \mathcal{K}(x^*) B^T \mathrm{Exp} \left(Z^T \mathrm{Ln}\left(\frac{x}{x^*}\right)\right), \quad \mathcal{K}(x^*) > 0
\end{equation}
This form will be the starting point for the analysis of balanced chemical reaction networks in the rest of this paper.

The matrix $B \mathcal{K}(x^*) B^T$ in (\ref{masterequation}) again defines a {\it weighted Laplacian} matrix for the complex graph, with weights given by the balanced reaction constants $\kappa_1(x^*), \cdots,$$ \kappa_r(x^*)$. Note that this is in general a {\it different} weighted Laplacian matrix than the one obtained before, cf. (\ref{Laplacian}). 
In particular, a main difference is that the weighted Laplacian $B \mathcal{K}(x^*) B^T$ is necessarily {\it symmetric}. Among others, cf. \cite{Bollobas}, this implies that the Laplacian $B \mathcal{K}(x^*) B^T$ is in fact independent of the {\it orientation} of the graph. Thus we may replace any reaction $ \mathcal{S} \rightleftharpoons \mathcal{P}$ by $\mathcal{P} \rightleftharpoons \mathcal{S}$ without altering the Laplacian $B \mathcal{K}(x^*) B^T$, in agreement with the usual understanding of a reversible reaction network.
The symmetrization of the Laplacian has been accomplished by the modification of $\mathrm{Ln}(x)$ into $\mathrm{Ln}\left(\frac{x}{x^*}\right)$, and using the assumption that $x^*$ is a thermodynamic equilibrium.

\smallskip

Note that $\mathcal{K}(x^*)$, and therefore the Laplacian matrix $B \mathcal{K}(x^*) B^T$, is in principle {\it dependent} on the choice of the thermodynamic equilibrium $x^*$. In the next Section \ref{subsec:linkage} we will see that actually this dependence is minor: for a connected complex graph the matrix $\mathcal{K}(x^*)$ is, {\it up to a positive multiplicative factor}, independent of the choice of the thermodynamic equilibrium $x^*$.

\subsection{The linkage classes of the complex graph}\label{subsec:linkage}
The complex graph provides a number of tools for the analysis of reaction networks. 
Recall that for any directed graph \cite{Bollobas}
\begin{equation}
\rank B = \rank L = c - \ell,
\end{equation}
where $c$ is the number of vertices of the graph, and $\ell$ is equal to the number of {\it components}\footnote{A directed graph is {\it connected} if there is a path (a number of un-oriented edges) between every two distinct vertices of the graph. The components of a directed graph are the maximal connected subgraphs.} of the complex graph
, the {\it linkage classes} in the terminology of \cite{HornJackson, Feinberg, Feinberg1}. 
Furthermore, if there is one linkage class in the network (i.e., the graph is connected and $\rank B = \rank L = c-1$), then 
\[
\ker L = \ker B^T = \spa \mathds{1}_c,
\]
where as before $\mathds{1}_c$ is the $c$-dimensional vector with elements all equal to $1$.

In general, if the reaction network has $\ell$ linkage classes then the network can be decomposed as follows. Let the $p$-th linkage class have $r_p$ reactions between $c_p$ complexes. Partition $Z$, $B$ and $\mathcal{K}$ matrices according to the various linkage classes as
\begin{equation}\label{decomposition}
\begin{array}{rcl}
Z &=& \begin{bmatrix}
Z_1 & Z_2 & \ldots & Z_{\ell}
\end{bmatrix}\\[4mm]
B &=& \begin{bmatrix}
B_1 & 0 & 0 & \ldots & 0\\
0 & B_2 & 0 & \ldots & 0\\
\vdots & \vdots & \ddots & \vdots & \vdots \\
0 & \ldots & 0 & B_{\ell -1} & 0\\
0 & \ldots & \ldots & 0 & B_{\ell}
\end{bmatrix}\\[16mm]
\mathcal{K}(x^*) &=& \text{diag}\big(\mathcal{K}_1(x^*),\mathcal{K}_2(x^*),\ldots,\mathcal{K}_{\ell}(x^*)\big)
\end{array}
\end{equation}
where for $p=1,\ldots,\ell$, $Z_p \in \bar{\mR}_+^{m \times c_p}$, $B_p \in \mR^{c_p \times r_p}$ and $\mathcal{K}_p \in \mR_+^{r_p \times r_p}$ denote, respectively, the complex stoichiometric matrix, incidence matrix and the diagonal matrix of balanced reaction constants for the $p$-th linkage class. Furthermore, $S_p=Z_pB_p$ is the stoichiometric matrix of the $p$-th linkage class. It follows that equation (\ref{masterequation}) can be expanded as
\begin{equation}\label{eq:linkage}
\dot{x}=-\sum_{p=1}^{\ell}Z_p B_p \mathcal{K}_p(x^*)B_p^T\Exp\left(Z_p^T\Ln\left(\frac{x}{x^*}\right)\right),
\end{equation}
expressing the contributions of each linkage class to the chemical reactions\footnote{Note that this yields a conserved moiety for the reaction network. In fact, the \emph{law of conservation of mass} states that for any complex stoichiometric matrix $Z$ there exists $u \in \mR_+^m$, such that $Z_p^{\top}u \in \spa \mathds{1}_{c_p}$ for $p=1,\ldots,\ell$. Since $B_p^T \mathds{1}_{c_p} =0$, this implies that $u^{\top}x$ is a conserved moiety for $\dot{x} = ZBv(x)$, for {\it all} forms of the reaction rate $v(x)$.}.

\smallskip

The above expressions also yield the following alternative characterization of the set of thermodynamic equilibria $\mathcal{E}$ given in (\ref{equilibria}). Since for a connected graph $\ker B^T = \spa \mathds{1}_c,$ the set $\mathcal{E}$ can be equivalently represented as, writing $S^T = B^TZ^T$,
\begin{equation}\label{equilibria1}
\mathcal{E}=\{x^{**} \in \mathbb{R}_+^m \mid Z^T \Ln \left(x^{**}\right) - Z^T \Ln \left(x^*\right) \in \spa \mathds{1}_c \}
\end{equation}
For non-connected graphs the same formula holds for every component.

This observation has the following useful consequence.
\begin{proposition}\label{other_kinetic_constants}
Consider the definition of the diagonal matrix $\mathcal{K}(x^*)$ of balanced reaction constants in (\ref{K}), and the formulation of the dynamics given in (\ref{masterequation}) for a particular thermodynamic equilibrium $x^*$. Assume that the complex graph is connected. Then for any other thermodynamic equilibrium $x^{**}$ there exists a positive constant $d^{**}$ such that
\[
\mathcal{K}(x^{**}) = d^{**}\mathcal{K}(x^*), \quad \Exp (Z^T \Ln \left(\frac{x}{x^{**}} \right)) = \frac{1}{d^{**}} \Exp(Z^T \Ln \left(\frac{x}{x^{*}} \right))
\]
Furthermore, for a non-connected complex graph there exist positive constants $d_p^{**}, p=1,\cdots, \ell,$ such that for each $p$-th linkage class the diagonal matrix $\mathcal{K}_p(x^*)$ given in (\ref{decomposition}) satisfies
\[
\mathcal{K}_p(x^{**}) = d_p^{**}\mathcal{K}_p(x^*), \quad p=1, \cdots, \ell
\]
\end{proposition}
\begin{proof}
Assume that the complex graph is connected. Then by (\ref{equilibria1}) there exists a constant $c^{**}$ such that $Z^T \Ln (x^{**}) = Z^T \Ln (x^*) + c^{**} \mathds{1}$. Hence by the definition of the balanced reaction constants $\kappa_j$ in (\ref{kappa}) we have 
\[
\begin{array}{rcl}
\kappa_j(x^{**}) & = & k_j^{\text{forw}} \exp \left( Z^T_{\mathcal{S}_j} \Ln(x^{**})\right) = k_j^{\text{forw}} \exp \left( 
Z^T_{\mathcal{S}_j} \Ln (x^{*}) + c^{**}\right) \\
& = & d^{**} k_j^{\text{forw}} \exp \left( Z^T_{\mathcal{S}_j} \Ln(x^{*})\right) = d^{**} \kappa_j(x^{*}),
\end{array}
\]
with $d^{**} := \exp (c^{**}) > 0$, for {\it each} $j=1, \cdots,r$. The rest follows easily.
\end{proof}
Hence the weighted Laplacian matrix $B \mathcal{K}(x^*) B^T$ is, up to a multiplicative factor, {\it independent} of the choice of the thermodynamic equilibrium $x^*$. The properties of this matrix (e.g., its eigenvalue distribution\footnote{Note that the set of eigenvalues of $B \mathcal{K}(x^*) B^T$ consists of $0$ (with multiplicity equal to the number of components) and strictly positive real numbers, cf. \cite{Bollobas}.}) therefore may serve as an important indicator of the network dynamics. In Section \ref{sec:modelreduction} we will exploit the Laplacian matrix for model reduction purposes.

\subsection{Deficiency}\label{subsec:deficiency}
An important notion to relate the structure of the complex graph to the {\it stoichiometry}, as introduced in the work of Feinberg \cite{Feinberg}, is the notion of {\it deficiency}. 
\begin{definition}\label{def:deficiency}
The deficiency of a chemical reaction network with complex stoichiometric matrix $Z$ and incidence matrix $B$ is defined as
\begin{equation}
\delta := \rank B - \rank{ZB} = \rank B - \rank S \geq 0
\end{equation}
A reaction network is said to have {\it zero-deficiency} if $\delta =0$. 
\end{definition}
\noindent
Note that zero-deficiency is equivalent to
\[
\ker Z \cap \im B =0,
\]
or to the mapping
$Z: \im B \subset \mathbb{R}^c \to \mathbb{R}^m$
being {\it injective}. Hence in the zero-deficiency case there is a one-to-one correspondence between the rate vector $\dot{x} \in \im S \subset \mathbb{R}^m$ of chemical species $x \in \mathbb{R}_+^m$ and the rate vector $\dot{y}\in \im B \subset \mathbb{R}^c$ of complexes $y \in \mathbb{R}^c$. Many chemical reaction networks are zero-deficient, although with growing complexity (especially in biochemical networks), deficiency greater than zero is likely to occur. This is also expressed in the following proposition showing that for a reaction network with zero-deficiency all its linkage classes also have zero-deficiency, but not necessarily the other way around.

\begin{proposition}\label{prop:zero_def_network}\cite{Othmer}
Consider a chemical reaction network with $\ell$ linkage classes as in (\ref{eq:linkage}). If the network has zero-deficiency then necessarily the linkage classes also have zero-deficiency. However, zero-deficiency of the linkage classes does not generally imply zero-deficiency of the total network. This {\it does} hold if and only if additionally
\[
\bigcap_{p=1}^{\ell} \im \bar{Z}_p \bar{B}_p= 0
\]
where $\bar{Z}_p$ and $\bar{B}_p$ denote the columns of $Z$, respectively of $B$, in (\ref{eq:linkage}) corresponding to $Z_p$, respectively $B_p$.
\end{proposition}

\subsection{Thermodynamical viewpoint}
In this section we will indicate the thermodynamical interpretation of the quantities introduced before, and show how this suggests a Lyapunov function which will be used in the next section. For more details regarding the thermodynamical approach to chemical reaction kinetics we refer to \cite{op,opk}.

Recall that for an ideal dilute solution the chemical potential $\mu_i$ of chemical species $i$ with mole number $n_i$ is given by
\begin{equation}
\mu_i(x_i) = \mu_i^o + RT \ln (\frac{n_i}{V}) = \mu_i^o + RT \ln (x_i)
\end{equation}
with $\mu_i^o$ a reference potential, $R$ a constant, $T$ the temperature, $V$ the volume, and $x_i = \frac{n_i}{V}$ the concentration. Equivalently, we have the inverse relation
\[
x_i = \exp \frac{(\mu_i - \mu_i^o)}{RT}
\]
The $m$-dimensional vector $\mu$ with components $\mu_i$ is called the {\it chemical potential} vector, while the vector $\mu^o$ with components $\mu_i^o$ is called the {\it reference chemical potential} vector.

Starting instead from the formulation of a balanced chemical reaction network in (\ref{masterequation}), corresponding to the thermodynamic equilibrium $x^*$, we may define the chemical potential vector $\mu$ and the reference chemical potential vector $\mu^0$ as
\begin{equation}
\mu(x) := RT \Ln \left(\frac{x}{x^*}\right), \quad \mu^o := - RT \Ln (x^*) 
\end{equation}
We conclude that $\mathrm{Ln} (\frac{x}{x^*}) = \mathrm{Ln}(x) - \mathrm{Ln}(x^*)$ is, up to the  constant $RT$, equal to the chemical potential vector, while $- \mathrm{Ln}(x^*)$ is, up to this same constant, equal to the reference chemical potential vector. 
Recalling from (\ref{wegscheider}) the expression of the vector of equilibrium constants as $K^{eq} = \Exp \left(S^T \Ln (x^*) \right)$, it follows that 
\begin{equation}
RT \Ln \left(K^{eq}\right) = (RT) S^T \Ln(x^*) = - S^T \mu^o
\end{equation}
An important role in equation (\ref{masterequation}) is formed by the quantity $\gamma(x) := Z^T \mu(x) =Z^T\Ln \left(\frac{x}{x^*}\right)$, which we will call the {\it complex thermodynamical affinity}. Correspondingly, we refer to $\gamma^o := Z^T \mu^o$ as the {\it reference complex thermodynamical affinity}\footnote{Note that $(RT) \Ln (K^{eq}) = (RT) S^T \Ln(x^*) = - B^T \gamma^o$, and thus the equilibrium constants correspond to {\it differences} of the reference complex thermodynamical affinities of the substrate and product complexes.}. The dynamics of a balanced reaction network is determined by this complex thermodynamical affinity $\gamma$, which acts as a `driving force' for the reactions. 
\begin{remark}
Consider on the other hand the vector $\alpha:= S^T \mu$, known as (minus\footnote{To remain in line with the usual sign convention in thermodynamics.}) the vector of {\it thermodynamical affinities}. Then for mass action kinetics It is not possible, see e.g. \cite{op,opk}, to express the vector of reaction rates $v(x)$ as a function of $\alpha$ (contrary to $\gamma$). This is illustrated in \cite{opk} by considering the reaction $X_1 \rightleftharpoons X_2$ with reaction rate $v(x_1,x_2) = k^{\text{forw}}x_1 - k^{\text{rev}}x_2$. When the concentrations $x_1$ and $x_2$ are doubled, then so is the reaction rate $v(x_1,x_2)$. However, the thermodynamical affinity $- \alpha$ remains the same.
\end{remark}
In thermodynamics the vector of chemical potentials is derived as the vector of partial derivatives of the {\it Gibbs' free energy}. This suggests to define the Gibbs' free energy as
\[
G(x) =RT x^T \mathrm{Ln}\left(\frac{x}{x^*}\right) + RT \left(x^* - x \right)^T \mathds{1}_m,
\]
where $\mathds{1}_m$ denotes a vector of dimension $m$ with all ones. Indeed, as is immediately checked
\[
\frac{\partial G}{\partial x}(x) = RT \mathrm{Ln}\left( \frac{x}{x^*}\right) = \mu(x)
\]
In the rest of this paper we will assume that $RT=1$, or equivalently, we will (re-)define
\begin{equation}\label{notation}
\begin{array}{l}
\mu(x) := \Ln \left( \frac{x}{x^*} \right),  \quad \gamma(x) := Z^T \mu(x) = Z^T \Ln \left( \frac{x}{x^*} \right) \\[3mm]
G(x):=x^T \mathrm{Ln}\left(\frac{x}{x^*}\right) + \left(x^{*} - x\right)^T \mathds{1}_m
\end{array}
\end{equation}
It follows that the equations of a balanced chemical reaction network (\ref{masterequation}) can be equivalently written as
\begin{equation}\label{masterequation1}
\dot{x} = - Z B \mathcal{K}(x^*) B^T \mathrm{Exp} \left(Z^T\frac{\partial G}{ \partial x}(x)\right), \quad \mu(x) = \frac{\partial G}{ \partial x}(x)
\end{equation}
In the next section, we will employ $G(x)$ as a {\it Lyapunov function}\footnote{Note that the definition of $G$ depends on the chosen thermodynamical equilibrium. Denoting the functions for different thermodynamic equilibria $x^*$ and $x^{**}$ by $G^*$, respectively $G^{**}$, it is seen that $G^{**}(x) = G^*(x) + x^T[ \Ln (x^*) - \Ln (x^{**})] + (x^{**} - x^*)^T \mathds{1}_m$.} for the chemical reaction network. 
In particular we will prove that $\frac{dG}{dt} \leq 0$; showing that $G$ is non-increasing along solution trajectories.

\medskip

A {\it geometrical interpretation} of Equations (\ref{masterequation}) and (\ref{masterequation1}) is as follows. Denote the dual space of the space of concentrations of chemical species $\mathcal{M} := \mathbb{R}^m$ by $\mathcal{M}^*$. Similarly, denote the dual space of $\mathcal{C} := \mathbb{R}^c$  by $\mathcal{C}^*$, and the dual of the space of reaction rates $\mathcal{R} = \mathbb{R}^r$ by $\mathcal{R}^*$. Define $v^*:=B^T \Exp(\gamma)$ and $y:=Bv(x)$. All ingredients of the equation (\ref{masterequation1}) are then summarized in the following diagram
\begin{equation}
\begin{array}{ccccccc}
& v \in \mathcal{R} & \overset{B}{\longrightarrow} & y \in \mathcal{C} & \overset{Z}{\longrightarrow} & x \in \mathcal{M}  &\\[6mm]
\mathcal{K}(x^*) & \mid &&& & \mid & G(x)  \\[6mm]
& v^* \in  \mathcal{R}^* & \underset{B^T}{\longleftarrow} & \gamma \in \mathcal{C}^* & \underset{Z^T}{\longleftarrow} & \mu \in \mathcal{M}^* & \\
& & & \circlearrowright &&&\\
&&& \mathrm{Exp} &&& \\
\end{array}
\end{equation}
which expresses the duality relations between all the variables involved. The concentration vector $x$ and its time-derivative $\dot{x}$ are elements of the linear space $\mathcal{M}$ with as conjugate vector the chemical potential vector $\mu \in \mathcal{M}^*$. They are related by the Gibbs' function $G(x)$ as $\mu = \frac{\partial G}{\partial x}(x)$. Furthermore, the vector $y$ is in the linear space $\mathcal{C}$, with conjugate vector the complex affinity $\gamma$. The relations between $y$ and $\dot{x}$ and dually between $\mu$ and $\gamma$ are given by $\dot{x}=Zy$, respectively $\gamma = Z^T \mu$. Also note that $\dot{x} = Zy = ZBv = Sv$, where the vector of fluxes $v$ is in the linear space $\mathcal{R}$, with conjugate vector $v^*:= - \big(\mathcal{K}(x^*)\big)^{-1}v \in \mathcal{R}^*$. The added complication in the diagram is the map $\mathrm{Exp}: \mathcal{C}^* \to \mathcal{C}^*$, which introduces a discrepancy between $v^*$ and $\alpha: = -B^T \gamma = -S^T \mu$.
%

\section{Equilibria and stability analysis of balanced reaction networks}\label{sec:ICSM}
In this section we show that the set of equilibria of a balanced chemical reaction network are actually equal to the set of thermodynamic equilibria $\mathcal{E}$, and we study their stability properties. 
\subsection{Equilibria of balanced reaction networks}
Making use of the formulation of the dynamics of balanced reaction networks in (\ref{masterequation}), we give a simple proof of the statement that all equilibria of a balanced reaction network are actually thermodynamic equilibria, and thus given by (\ref{equilibria}). A similar result was obtained in the classical papers \cite{Horn, HornJackson, Feinberg1} for a different class of chemical reaction networks (roughly speaking, weakly reversible  networks of deficiency zero or deficiency one under additional conditions).
\begin{theorem}\label{th:char}
Consider a balanced chemical reaction network $\dot{x} = Sv = ZBv$ with $m$ species and $r$ reactions governed by mass action kinetics, with thermodynamic equilibrium $x^*$, i.e., $v(x^*)=0$, described as in (\ref{masterequation}).
Then the set of all equilibria is equal to the set $\mathcal{E}:=\{x^{**} \in \mathbb{R}_+^m \mid S^T \Ln (x^{**})= S^T \Ln (x^*)=0\}$ of thermodynamic equilibria given in (\ref{equilibria}).
\end{theorem}
\begin{proof}
%
%
%
%
%
Denote for the $j$-th reaction the substrate complex by $\mathcal{S}_j$ and the product complex by $\mathcal{P}_j$. Let as before $Z_{\mathcal{S}_j}$ and $Z_{\mathcal{P}_j}$ denote the columns of the complex stoichiometric matrix $Z$ of the reaction network, corresponding respectively to substrate complex $\mathcal{S}_j$ and product complex $\mathcal{P}_j$. Define as before
\[
\mu(x)=\mathrm{Ln}\left( \frac{x}{x^*}\right), \, \gamma(x):=Z^T \mu(x), \, \gamma_{\mathcal{S}_j}(x) =Z^T_{\mathcal{S}_j} \mu(x), \,\gamma_{\mathcal{P}_j}(x)=Z^T_{\mathcal{P}_j} \mu(x) 
\]
Suppose $x^{**}$ is an equilibrium, i.e.,
$ZB\mathcal{K}(x^*) B^T \mathrm{Exp} \left(Z^T \mathrm{Ln}\left(\frac{x^{**}}{x^*}\right)\right) = 0,$
or equivalently $ZB\mathcal{K}(x^*) B^T \mathrm{Exp} \left(Z^T \mu(x^{**})\right) = 0$. Then also
\[
\mu^T(x^{**})ZB\mathcal{K}(x^*) B^T \mathrm{Exp} \left(Z^T \mu(x^{**})\right) =0 
\]
Denoting the columns of $B$ by $b_1, \cdots, b_r,$ we have
\[
B\mathcal{K}(x^*) B^T = \sum_{j=1}^r \kappa_j(x^*)b_jb_j^T
\]
while 
\[
\begin{array}{l}
\mu^T(x^{**})Zb_j = \mu^T(x^{**}) \left(Z_{\mathcal{P}_j} - Z_{\mathcal{S}_j}\right) = \gamma^T_{\mathcal{P}_j}(x^{**}) -  \gamma^T_{\mathcal{S}_j}(x^{**}), \\[4mm]
b_j^T \Exp\left(Z^T \mu(x^{**}) \right) = \exp \left(\gamma^T_{\mathcal{P}_j}(x^{**})\right) - \exp \left(\gamma^T_{\mathcal{S}_j}(x^{**})\right) 
\end{array}
\]
It follows that
\[
\begin{array}{rl}
0 & = \mu^T(x^{**})ZB\mathcal{K}(x^*) B^T \mathrm{Exp} \left(Z^T \mu(x^{**})\right) \\[2mm]
&= \sum_{j=1}^r
\left[\gamma_{\mathcal{P}_j}(x^{**}) - \gamma_{\mathcal{S}_j}(x^{**})\right] \left[\exp\big(\gamma_{\mathcal{P}_j}(x^{**}) \big)-\exp\big(\gamma_{\mathcal{S}_j}(x^{**}) \big)\right]\kappa_j(x^*)
\end{array}
\]
Since the exponential function is a strictly increasing function and $\kappa_j(x^*) >0, j=1, \cdots, r,$ this implies that all terms in the summation are zero and thus
\[
\begin{array}{ll}
\gamma_{\mathcal{P}_j}(x^{**}) = \gamma_{\mathcal{S}_j}(x^{**}),   & \quad j=1, \cdots, r
\\[2mm]
\exp\big(\gamma_{\mathcal{P}_j}(x^{**}) \big) = \exp\big(\gamma_{\mathcal{S}_j}(x^{**}) \big),  &\quad j=1, \cdots, r
\end{array}
\]
which is equivalent to
\[
B^T \gamma(x^{**}) = B^TZ^T \mu(x^{**}) = 0, \quad B^T \Exp \left(\gamma(x^{**})\right) = B^T \Exp \left(Z^T \mu(x^{**})\right) = 0
\]
The first equality tells us that $x^{**} \in \mathcal{E}$, and thus a thermodynamical equilibrium (as also follows from the second equality).
\end{proof}  


\subsection{Asymptotic stability}
In the next theorem we show that $G$ serves as a Lyapunov function for (\ref{masterequation}), implying global asymptotic stability with respect to the set of equilibria $\mathcal{E}$ given by (\ref{equilibria}), i.e., $\mathcal{E}$ is globally attractive. 

\begin{theorem}\label{tm:Lyap}
Consider a balanced mass action reaction network given by (\ref{masterequation}) or, equivalently, by (\ref{masterequation1}). Then for every initial condition $x(0)\in \mR_+^{m}$, the species concentration $x$ converges for $t \to \infty$ to $\mathcal{E}$.
\end{theorem}

\begin{proof}
First we prove that the function $G$ in (\ref{notation}) satisfies
\begin{equation}\label{lyap}
G(x^{*}) =0, \quad G(x) > 0 , \quad \forall x \neq x^{*},
\end{equation}
and is {\it proper}, i.e., for every real $ c >0$ the set $\{x\in \bar{\mR}_+^{m} \mid G(x) \leq c\}$ is compact. 
With regard to (\ref{lyap}) we first note the following. 
Let $x_i$ and $x_i^*$ denote the $i$-th elements of $x$ and $x^*$ respectively. From the strict concavity of the logarithmic function 
\begin{equation}\label{eq:z}
z- 1\geq \ln(z), \quad \forall z \in \mR_+, 
\end{equation}
with equality if and only if $z=1$.
Putting $z=\frac{x_i^*}{x_i}$ in equation (\ref{eq:z}), we get
\[
x_i^*-x_i+x_i \ln \left(\frac{x_i}{x_i^*}\right) \geq 0,
\] 
with equality if and only if $x_i = x_i^*$. This implies that
\[
G(x)=\sum_{i=1}^m\left[x_i^*-x_i+x_i\ln\left(\frac{x_i}{x_i^*}\right)\right] \geq 0,
\]
with equality if and only if $x_i = x_i^*, i=1, \cdots, m$. Thus $G$ has a strict minimum at $x=x^*$ and satisfies (\ref{lyap}). Properness of $G$ is readily checked.

Next we will show that $\dot{G}(x):= \frac{\partial^T G}{\partial x}(x) Sv(x) = \frac{dG}{dt}(x)$ satisfies
\begin{equation}\label{lyap1}
\dot{G}(x) \leq 0 \qquad \forall x \in \mR_+^{m},
\end{equation}
and
\begin{equation}\label{lyap2}
\dot{G}(x) = 0 \mbox{ if and only if } x\in \mathcal{E}.
\end{equation}
As in the proof of Theorem \ref{th:char}, let the $j$-th reaction be between substrate complex $\mathcal{S}_j$ and product complex $\mathcal{P}_j$, and let $Z_{\mathcal{S}_j}$ and $Z_{\mathcal{P}_j}$ denote the columns of the complex stoichiometric matrix $Z$ corresponding to complexes $\mathcal{S}_j$ and $\mathcal{P}_j$. Define as before
\[
\mu(x)=\mathrm{Ln}\left( \frac{x}{x^*}\right), \, \gamma(x):=Z^T \mu(x), \, \gamma_{\mathcal{S}_j}(x)=Z^T_{\mathcal{S}_j} \mu(x), \,\gamma_{\mathcal{P}_j}(x) =Z^T_{\mathcal{P}_j} \mu(x)
\] 
Then compute\footnote{See the proof of Theorem \ref{th:char} for an analogous reasoning.}
\begin{equation}\label{sum}
\begin{array}{rl}
\dot{G}(x) &= \frac{\partial^T G}{\partial x}(x)\dot{x}= -\mu^T(x) ZB\mathcal{K}(x^*)B^T \Exp(Z^T\mu(x)) \\[2mm]
&= -\gamma^T(x) B\mathcal{K}(x^*)B^T \Exp(\gamma(x)) \\[2mm] 
&= \sum_{j=1}^r
\left[\gamma_{\mathcal{P}_j}(x) - \gamma_{\mathcal{S}_j}(x)\right] \left[\exp\big(\gamma_{\mathcal{P}_j}(x) \big)-\exp\big(\gamma_{\mathcal{S}_j}(x) \big)\right]\kappa_j(x^*)\\[2mm]
& \leq 0 \, ,
\end{array}
\end{equation}
since $\kappa_j(x^*) >0$ for $j=1,\ldots, r,$ and the exponential function is strictly increasing. 
The ``if'' part in (\ref{lyap2}) is trivial. For the ``only if'' part, the summand in the third line of (\ref{sum}) is zero only if $\gamma_{\mathcal{S}_j}(x) - \gamma_{\mathcal{P}_j}(x) =0$ for every $j$. This is equivalent to having $B^T\gamma(x) = 0$. Thus, $\dot{G}(x)=0$ only if $B^T\gamma=B^TZ^T\Ln\left(\frac{x}{x^*}\right)=0$. It follows that $\dot{G}(x) = 0$ if and only if $x\in \mathcal{E}$.

Suppose that the set $\mR_+^m$ is {\it not} forward invariant with respect to (\ref{masterequation1}), and that $x_{i}(t)=0$ for some $t$ and $i \in \{1, \cdots, m\}$. Without loss of generality assume that the species with concentration $x_i$ is taking part in at least one complex which is involved in a reaction. Let $\mathcal C_i$ denote the subset of complexes which contain $x_i$, i.e., $Z_{i,\mathcal{C}_i}\neq 0$. Then it follows that $\Exp(Z_{\mathcal{C}_i}^T\mu(x(t))) = 0$ and thus 
\begin{equation*}
\dot x_{i}(t)  = - Z^{i} B \mathcal{K}(x^*) B^T \Exp(Z^T\mu(x)) > 0,
\end{equation*}  
where $Z^i$ is the $i$-th row vector of $Z$. This inequality is due to the fact that the terms corresponding to the positive $i$-th diagonal element of the weighted Laplacian matrix $B \mathcal{K}(x^*) B^T$ are all zero, while there is at least one term corresponding to a non-zero, and therefore strictly negative, off-diagonal element of $B \mathcal{K}(x^*) B^T$. This implies\footnote{A similar argument can be used directly for the mass action kinetics rate equations.} that $\mR_+^m$ must be forward invariant. 

Since $G$ is proper (in $\bar{\mR}_+^{m}$) and the state trajectory $x(\cdot)$ remains in $\mR_+^m$, (\ref{lyap1}) implies that $x(\cdot)$ is bounded in $\mR_+^m$. Therefore, boundedness of $x(\cdot)$, together with equations (\ref{lyap1}) and (\ref{lyap2}), imply that the species concentration $x$ converges to $\mathcal{E}$ by an application of LaSalle's invariance principle. 
\end{proof}
\begin{remark} A similar reasoning for showing that $\dot{G} \leq 0$ was used in \cite{AngeliEJC}, see also \cite{Horn, Feinberg, Feinberg1}; however making use of the {\it convexity} of the exponential function instead of the weaker property that the exponential function is {\it increasing} as in our proof. In this respect, it can be noted that Theorem \ref{tm:Lyap} remains unaltered if we replace the dynamics (\ref{masterequation}) by any equations of the form (not corresponding anymore to mass action kinetics)
\[
\dot{x} = - Z B \mathcal{K}(x^*) B^T F \left(Z^T \mathrm{Ln}\left(\frac{x}{x^*}\right)\right),
\]
with $F: \mathbb{R}^c \to \mathbb{R}^c$ a mapping $F(y_1,\cdots,y_c)= \diag (f_1(y_1), \cdots,f_c(y_c))$, where the functions $f_i, i=1,\cdots,c,$ are all monotonously increasing.
\end{remark}

\subsection{Equilibrium concentration corresponding to an initial concentration}
In this section, we show that for every initial concentration vector $x_0 \in \mR_+^m$, the solution trajectory of (\ref{masterequation}) converges to a {\it unique} thermodynamic equilibrium in $\mathcal{E}$.
Our reasoning is very similar to the proof of zero-deficiency theorem provided in \cite{Feinberg1} and is based on the following proposition in there. Recall from the Introduction that the product $x \cdot z \in \mathbb{R}^m$ is defined as the element-wise product $(x \cdot z)_i= x_i z_i, \, i=1, \cdots,m$.

\begin{proposition}\label{prop:uniq}
Let $U$ be a linear subspace of $\mR^m$, and let $x^*,x_0 \in \mR_+^m$. Then there is a unique element $\mu \in U^{\perp}$, such that $\big(x^* \cdot \Exp(\mu)-x_0\big) \in U$.
\end{proposition}
\begin{proof}
See proof of \cite[Proposition B.1, pp. 361-363]{Feinberg1}.
\end{proof} 
\begin{theorem}\label{th:unique}
Consider the balanced chemical reaction network (\ref{masterequation}). Then for every $x_0 \in \mR_+^m$, there exists a unique $x_1 \in \mathcal{E}$ such that the solution trajectory of (\ref{masterequation}) starting from $x_0$ converges for $t \to \infty$ to $x_1$. Hence there exists a surjective map $\chi:\mR_+^m \to \mathcal{E}$ that assigns to every initial state its asymptotic thermodynamic equilibrium\footnote{The map $\chi$ is an analytic function, following a similar argument as used in \cite[Theorem 6]{Sontag} which deals with weakly-reversible zero-deficiency chemical networks.}.
\end{theorem}
\begin{proof}
With reference to Proposition \ref{prop:uniq}, define $U=\spa S$, and observe that $U^{\perp}=\ker S^T$. Let $x^*,x_0\in \mR_+^m$. Then by Proposition \ref{prop:uniq}, there exists a unique $\mu \in \ker S^T$ such that $x^* \cdot \Exp(\mu)-x_0 \in \spa S$. Define $x_1:=x^* \cdot \Exp(\mu) \in \mR_+^m$. It follows that $S^{T}\mu=S^T\Ln\left(\frac{x_1}{x^*}\right)=0$, i.e., $x_1 \in \mathcal{E}$. Also, since $x_1-x_0 \in \spa S$, $x_1 \in \{\xi\in\mR^m \mid \xi - x_0 \in \spa S\}$ which is an invariant set of the dynamics.
Together with Theorem \ref{tm:Lyap}, the state trajectory with initial condition $x_0$ converges to the unique equilibrium $x_1 \in \mathcal{E}$. 
\end{proof}

\subsection{Relation with consensus dynamics}
The form of the equations (\ref{sontag}) is reminiscent of {\it consensus dynamics} for multi-agent systems with non-symmetric communication. In fact, a standard consensus algorithm in the vector $y \in \mathbb{R}^c$ can be written as \cite{Murray}
\[
\dot{y} = -L_{\text{cons}} y,
\]
where $L_{\text{cons}}$ is the weighted Laplacian matrix of the communication graph. This weighted Laplacian is constructed similarly as in (\ref{Laplacian}), with the difference that $\Delta$ is now defined as the diagonal matrix whose $(\rho,\rho)$-th element is equal to the sum of the elements of the $\rho$-th {\it row} of $A$, implying that the sum of the elements of every {\it row} of $L_{\text{cons}}$ is zero. Note that in contrast for the weighted Laplacian $L$ in (\ref{Laplacian}) the sum of the elements of every {\it column} is zero\footnote{Of course, another main difference between the two dynamics is the appearance of the complex stoichiometric matrix $Z$ on the left of the right-hand side of the differential equation (\ref{sontag}), and the much more involved term $\mathrm{Exp} \left(Z^T \mathrm{Ln}(x)\right)$ instead of simply $y$ (see however \cite{Cortes}).}.
Seen from this perspective, the equations (\ref{masterequation}) for balanced chemical reaction networks, exhibiting the {\it symmetric} weighted Laplacian matrix  $B\mathcal{K}(x^*)B^T$, correspond to the consensus algorithm for multi-agent systems with {\it symmetric} communication structure \cite{Murray}.

Finally, Theorem \ref{th:unique} is closely related to the property that the asymptotic consensus value in consensus dynamics is a function of the initial state. In fact, the surjective map $\chi:\mR_+^m \to \mathcal{E}$ which assigns to every initial state its asymptotic equilibrium is similar to the $\chi$-function in the $\chi$-consensus algorithm of \cite{Cortes}. 

\section{Chemical reaction networks with boundary fluxes and interconnection of reaction networks}\label{sec:BF}
As discussed in Section 2.1, the interaction of a chemical reaction network with the environment or another chemical reaction network can be modeled by identifying a vector of {\it boundary chemical species} $x_b \in \mathbb{R}^b$, which is a sub-vector of the vector $x$ of concentrations of all the chemical species in the network. These boundary chemical species are the species that are subject to inflow or outflow boundary fluxes $v_b$. The natural conjugate vector is (up to the constant $RT$) the vector of chemical potentials $\mu_b=S_b^T\mu \in \mathbb{R}^b$ of {\it boundary potentials}. Indeed, up to the factor $RT$, the product $\mu_b^Tv_b$ is equal to the {\it power} entering or leaving the chemical reaction network due to the influx or efflux of boundary chemical species. 

By combining (\ref{stoichiometry}) and (\ref{masterequation}) this leads to the following equations
\begin{equation}\label{io}
\begin{array}{rcl}
\dot{x} & = &   - Z B \mathcal{K}(x^*) B^T \mathrm{Exp} \left(Z^T \Ln\left(\frac{x}{x^*}\right)\right) + S_bv_b, \\[2mm] 
\mu_b & = & S_b^T \Ln\left(\frac{x}{x^*}\right) ,
\end{array}
\end{equation}
which define an input-state-output system with inputs\footnote{Note however that the elements of $v_b$ can be either {\it in}flows or {\it out}flows for the network.} $v_b$ and outputs $\mu_b$. Let $G$ be defined by equation (\ref{notation}). Define as before $\gamma(x):=Z^{\top}\Ln\left(\frac{x}{x^*}\right)$. As an immediate consequence we obtain the following energy balance
\begin{equation}
\frac{d}{dt} G(x) = -\gamma^T(x) B \mathcal{K}(x^*) B^T \mathrm{Exp}(\gamma(x)) + \mu_b^Tv_b 
\end{equation}
Furthermore, since the exponential function is increasing it is easily verified along the lines of the proof of Theorem \ref{tm:Lyap} \big(see in particular (\ref{sum})\big), that
$
\gamma^T \mathcal{K} B^T \mathrm{Exp}(\gamma) \geq 0,
$
thus showing the {\it passivity} property
\begin{equation}
\frac{d}{dt} G \leq \mu_b^Tv_b 
\end{equation}

\subsection{Interconnection of chemical reaction networks}
The most basic way of interconnecting chemical reaction networks is through {\it shared boundary chemical species}. Consider for simplicity of notation the interconnection of {\it two} chemical reaction networks which have all their boundary chemical species in common; the generalization to the general case is straightforward. Let $B_i$ denote the incidence matrix of the complex graph, and $Z_i$ the complex stoichiometry matrix of network $i=1, 2$. The complex graph of the interconnected chemical reaction network is the direct union of the complex graphs of networks $1$ and $2$, with incidence matrix $B$ given as the direct product
\begin{equation}\label{B_interconnect}
B = \begin{bmatrix} B_1 & 0 \\ 0 & B_2 \end{bmatrix}
\end{equation}
Note that the complex graph of such an interconnected network is {\it not} connected. In fact, if all the constituent networks are connected then their complex graphs are {\it linkage classes} of the interconnected network, cf. Section \ref{subsec:deficiency}. 

The complex stoichiometry matrix of the interconnected network is more involved. Permute the chemical species $x_1, x_2$ such that  
\[
x_1=(\hat{x}_{1},x_{b1}), x_2 = (x_{b2},\hat{x}_{2}), \quad Z_1 = \begin{bmatrix} \hat{Z}_1 \\ Z_{b1} \end{bmatrix}, \quad Z_2 = \begin{bmatrix} Z_{b2} \\ \hat{Z}_2  \end{bmatrix},
\]
where $Z_{b1}, Z_{b2}$ are matrices with the same number of rows, equal to the number of shared boundary species 
\[
x_b:=x_{b1} = x_{b2} \in \mathbb{R}^b
\]
Then the stoichiometry matrix $Z$ of the interconnected network is given as
\begin{equation}\label{stoichint}
Z = \begin{bmatrix} \hat{Z}_1 & 0 \\
Z_{b1} & Z_{b2} \\ 0 & \hat{Z}_2 \end{bmatrix}
\end{equation}
Note that in general the property of {\it zero-deficiency} is not maintained under interconnection. In fact, this was already discussed in Proposition \ref{prop:zero_def_network} in Section \ref{subsec:deficiency} following \cite{Othmer}.
\begin{remark}\label{remarkinterconnection}
It may occur that one or more columns of the matrix $Z$ in (\ref{stoichint}) are actually the same, corresponding to the case that the two reaction networks have {\it shared complexes} (by definition only consisting of shared boundary chemical species; i.e., for which the corresponding columns of $ \hat{Z}_1$ and $ \hat{Z}_2$ are zero). Then one may {\it identify} the equal columns in the matrix $Z$, and thus obtain a reduced network with complex graph consisting of the union of the complex graphs of the two networks with the vertices corresponding to the shared complexes being identified. In this case the complex graph of the interconnected network may still be connected.
\end{remark}

\smallskip

The property of 'balancedness' is not necessarily maintained under interconnection. Sufficient conditions for this to happen are given in the following proposition.
\begin{proposition}\label{prop:int}
The interconnection of two balanced reaction networks with vectors of equilibrium constants $K^{eq}_1, K^{eq}_2$ is again balanced if and only if there exists $(\hat{x}_1^{**}, \hat{x}_2^{**}, x_b^{**})$ such that
\begin{equation}\label{interconnectedequilibrium}
\Ln\left(\begin{bmatrix} K^{eq}_1 \\ K^{eq}_2 \end{bmatrix}\right) = \begin{bmatrix} B_1^T\hat{Z}_1^T & 0 & B_1^TZ_{b1}^T \\ 0 & B_2^T\hat{Z}_2^T & B_2^TZ_{b2}^T \end{bmatrix}\Ln\left(\begin{bmatrix} x_1^{**} \\ x_2^{**} \\ x_b^{**} \end{bmatrix}\right)
\end{equation}
There exists such a thermodynamic equilibrium for the interconnected network if there is a partition $\{1,\cdots,b\} = I_1 \bigcup I_2$ such that all columns of $B_1^TZ_{b1}^T$ corresponding to the index set $I_1$ are contained in the image of $B_1^T\hat{Z}_{1}^T$, and all columns of $B_2^TZ_{b2}^T$ corresponding to the complementary index set $I_2$ are contained in the image of $B_2^T\hat{Z}_{2}^T$. 
\end{proposition}

\begin{proof}
By assumption there exist thermodynamic equilibrium points $(\hat{x}_1^{*}, x_{b1}^{*})$ and $(x_{b2}^*, \hat{x}_2^{*})$ of the two individual networks. That means (cf. \ref{wegscheider})) that
\begin{align}
\label{Keq1} \Ln(K^{eq}_1) & = B_1\begin{bmatrix} \hat{Z}_1^T & Z_{b1}^T \end{bmatrix}\Ln \left( \begin{bmatrix} \hat{x}_1^{*} \\ x_{b1}^{*} \end{bmatrix} \right) \\
\label{Keq2} \Ln(K^{eq}_2) & = B_2\begin{bmatrix} \hat{Z}_2^T & Z_{b2}^T \end{bmatrix}\Ln \left( \begin{bmatrix} \hat{x}_2^{*} \\ x_{b2}^{*} \end{bmatrix} \right)
\end{align} 
Now define $x_b^{**} \in \mathbb{R}^b$ by taking its $i$-th component for $i \in I_1$ equal to the $i$-th component of $x_{b2}^*$, and for $i \in I_2$ equal to the $i$-th component of $x_{b1}^*$. Then there exist corresponding $\hat{x}_1^{**}, \hat{x}_2^{**}$ such that (\ref{interconnectedequilibrium}) is satisfied.
\end{proof}
Under the assumptions of Proposition \ref{prop:int}) it follows that the balanced interconnected network is given as
\begin{equation}
\begin{bmatrix} \dot{\hat{x}}_1 \\ \dot{\hat{x}}_2 \\ \dot{x}_{b} \end{bmatrix} = - ZB\begin{bmatrix} \mathcal{K}_1(\hat{x}_1^{**},x_b^{**}) & 0 \\ 0 & \mathcal{K}_2(\hat{x}_2^{**},x_b^{**})\end{bmatrix} B^T \Exp\left( Z^T \Ln\left(\begin{bmatrix} \hat{x}_1/{\hat{x}}_1^{**} \\ \hat{x}_2/\hat{x}_2^{**} \\ x_b/{x}_{b}^{**} \end{bmatrix}\right) \, \right), 
\end{equation}
where $B$ as in (\ref{B_interconnect}), $Z$ as in (\ref{stoichint}), $\hat{x}_1^{**}, \hat{x}_2^{**},x_b^{**}$ satisfy (\ref{interconnectedequilibrium}), and $\mathcal{K}_1(\hat{x}_1^{**},x_b^{**})$ and $\mathcal{K}_2(\hat{x}_2^{**},x_b^{**})$ are proportionally related to $\mathcal{K}_1(\hat{x}_1^{*},x_{1b}^{*})$ and $\mathcal{K}_2(\hat{x}_2^{*},x_{2b}^{*})$, respectively, as follows from Proposition \ref{other_kinetic_constants}.

\subsection{Interconnection arising from port interconnection}
The above procedure for interconnection of chemical reaction networks can be also interpreted as arising from {\it power-port interconnection} at the boundary chemical species. Recall the formulation (\ref{io}) of an open chemical reaction network with inputs $v_b$ being the influx/efflux of the boundary chemical species and $\mu_b$ their chemical potentials. Then the interconnection of two chemical reaction networks (as above under the simplifying assumption that all boundary chemical species are shared) can be seen to result from the power-port interconnection constraints
\begin{equation}\label{interconnection}
\mu_{b1} = \mu_{b2}, \quad v_{b1} + v_{b2} =0,
\end{equation}
expressing that the chemical potentials of the boundary chemical species are equal, while the boundary fluxes of the two networks add up to zero (conservation of boundary chemical species).
Indeed, the resulting interconnected dynamics given by the differential-algebraic equations
\[
\begin{array}{rcl}
\begin{bmatrix} \dot{x}_1 \\[2mm] \dot{x}_2 \end{bmatrix} & = &
- \begin{bmatrix} -Z_1B_1\mathcal{K}_1(x_1^*)B_1^T \Exp (Z_1^T \mu_1(x_1)) \\[2mm]
Z_2B_2\mathcal{K}_2(x_2^*)B_2^T \Exp (Z_2^T \mu_2(x_2)) \end{bmatrix} + \begin{bmatrix} I \\[2mm] -I \end{bmatrix}v \\[5mm]
0 & = & \begin{bmatrix} I & -I \end{bmatrix} \begin{bmatrix} \mu_{b1}(x_1) \\[2mm] \mu_{b2}(x_2) \end{bmatrix}
\end{array}
\]
with $v = v_{b1} = - v_{b2}$ acting as a vector of Lagrange multipliers, can be seen, after elimination of the algebraic constraint $\mu_{b1}(x_1) = \mu_{b2}(x_2)$ and the Lagrangian multipliers $v$, to result in the dynamics
\[
\dot{x} =  - ZB\mathcal{K}(x_1^*,x_2^*)B^T \Exp (Z^T \mu(x_1,x_2))
\]
where $x = (\hat{x}_1, \hat{x}_2, x_b)$ denotes the total vector of chemical species (internal species of both networks, together with the shared boundary species).
Here the complex stoichiometric matrix $Z$ and incidence matrix $B$ of the interconnected network are as above, and where $\mathcal{K}(x_1^*,x_2^*) = \mathrm{diag} (\mathcal{K}_1(x_1^*), \mathcal{K}_2(x_2^*))$ and $\mu(x_1,x_2) = \begin{bmatrix} \mu_1(x_1) \\ \mu_2(x_2) \end{bmatrix}$. 
Recalling that $\mu_b^Tv_b$ is (modulo the constant $RT$) the {\it power} provided to the chemical reaction network, this implies that the interconnection (\ref{interconnection}) is {\it power-preserving}, that is $\mu_{b1}^Tv_{b1} + \mu_{b2}^Tv_{b2} =0$, in line with the standard way of interconnecting physical networks \cite{op, opk, opcircuits}.

\section{Model reduction of chemical reaction networks}\label{sec:modelreduction}
For many purposes one may wish to {\it reduce} the number of dynamical equations of a complex chemical reaction network in such a way that the behavior of a number of key chemical species is approximated in a satisfactory way. 

Furthermore, it is important that the reduced system again allows the interpretation of a chemical reaction network (but with much less chemical species involved). In fact, most of the approaches to model reduction of large-scale chemical reaction networks aim at doing exactly this: they simplify the pathways of the chemical reaction network by leaving out intermediate species or complexes; see e.g. \cite{Farkas, Anderson}, or the standard reduction of mass action enzymatic reactions to their Michaelis-Menten description \cite{klipp}.

In the following we will propose a model reduction approach which is based on a {\it reduction of the complex graph}. First we recall from \cite{vdsSCL} the following result regarding Schur complements of weighted Laplacian matrices.
\begin{proposition}
Consider a directed graph with vertex set $\mathcal{V}$ and with weighted Laplacian matrix $L=B\mathcal{K}(x^*)B^T$. Then the Schur complement with respect to any subset of vertices $\mathcal{V}_r \subset \mathcal{V}$ is well-defined, and is the weighted Laplacian matrix $\hat{L} = \hat{B} \hat{\mathcal{K}}(x^*) \hat{B}^T$ of a reduced directed graph with incidence matrix $\hat{B}$ whose vertex set is equal to $\mathcal{V} - \mathcal{V}_r$. 
\end{proposition}
\begin{remark}
Note that the reduced graph defined by $\hat{B}$ may contain additional edges between the remaining vertices in  $\mathcal{V} - \mathcal{V}_r$. Furthermore, the weights of the reduced diagonal weight matrix $\hat{\mathcal{K}}(x^*)$ are in general {\it different} from the weights of the original weight matrix $\mathcal{K}(x^*)$, and represent {\it equivalent} balanced reaction constants\footnote{In electrical circuit theory the process of reduction of a resistive network to a resistive network with fewer vertices but equivalence resistance is sometimes referred to as {\it Kron reduction}, cf. \cite{kron, Dorfler}.}.
\end{remark}
This result can be directly applied to the complex graph, yielding a reduction of the chemical reaction network by reducing the number of complexes. The procedure is as follows. Consider a balanced reaction network described in the standard form (\ref{masterequation})
\[
\Sigma: \quad \dot{x} = - Z B \mathcal{K}(x^*) B^T \mathrm{Exp} \left(Z^T \mathrm{Ln}\left(\frac{x}{x^*}\right)\right)
\]
Reduction will be performed by {\it deleting certain complexes in the complex graph}, resulting in a reduced complex graph with weighted Laplacian $\hat{L} = \hat{B} \hat{\mathcal{K}}(x^*) \hat{B}^T$. Furthermore, leaving out the corresponding columns of the complex stoichiometric matrix $Z$ one obtains a reduced complex stoichiometric matrix $\hat{Z}$ (with as many columns as the remaining number of complexes in the complex graph), leading to the reduced reaction network
\begin{equation}\label{reduced}
\hat{\Sigma}: \quad \dot{x} = - \hat{Z} \hat{B} \hat{\mathcal{K}}(x^*) \hat{B}^T \mathrm{Exp} \left(\hat{Z}^T \mathrm{Ln}\left(\frac{x}{x^*}\right)\right), \quad \hat{\mathcal{K}}(x^*) > 0
\end{equation}
Note that $\hat{\Sigma}$ is again a {\it balanced chemical reaction network} governed by mass action kinetics, with a reduced number of complexes and with, in general, a different set of reactions (edges of the complex graph). Note furthermore that the thermodynamic equilibrium $x^*$ of the original reaction network $\Sigma$ is a thermodynamic equilibrium of the reduced network $\hat{\Sigma}$ as well. 

Furthermore, we can derive the following properties relating $\Sigma$ and $\hat{\Sigma}$.

\begin{proposition}\label{prop:reduction}
Consider the balanced reaction network $\Sigma$ and its reduced order model $\hat{\Sigma}$ given by (\ref{reduced}). Denote their sets of equilibria by $\mathcal{E}$, respectively $\hat{\mathcal{E}}$. Then $\mathcal{E} \subset \hat{\mathcal{E}}$. Furthermore, if $\Sigma$ has deficiency zero then so does $\hat{\Sigma}$.
\end{proposition}
\begin{proof}
Assume that the complex graph is connected; otherwise the same argument can be repeated for every component (linkage class). 
Recall from (\ref{equilibria}) that the set of equilibria $\mathcal{E}$ is given as
$\{x^{**} \mid \Ln^T(x^{**})S = \Ln^T(x^*)S \},$ 
where $S=ZB$. It follows that $\mathcal{E}$ is equivalently represented as
\[
\mathcal{E} = \{x^{**} \mid \Ln^T (x^{**})ZL = \Ln^T(x^*)ZL \},
\]
where $L:=B \mathcal{K}(x^*) B^T$ is the weighted Laplacian matrix of the complex graph. Now consider a subset $\mathcal{V}_r$ of the set of all complexes, which we wish to leave out in the reduced network. After permutation of the complexes we partition $L$ and $Z$ correspondingly as
\[
L = \begin{bmatrix} L_{11} & L_{12} \\ L_{21} & L_{22} \end{bmatrix}, \quad Z = \begin{bmatrix} Z_1 & Z_2 \end{bmatrix}
\]
where $\mathcal{V}_r$ corresponds to the last part of the vertices. Then post-multiply the equations
$\Ln^T (x^{**})ZL = \Ln^T(x^*)ZL$ by the invertible matrix
\begin{equation}\label{post}
\begin{bmatrix} I & 0 \\ -L_{22}^{-1}L_{21} & I \end{bmatrix} ,
\end{equation}
not changing the solution set $\mathcal{E}$. It follows that
\[
\begin{array}{rcl}
\mathcal{E} &=  \{x^{**} \mid & \Ln^T (x^{**}) \begin{bmatrix} Z_1 & Z_2 \end{bmatrix} 
\begin{bmatrix} L_{11} - L_{12}L_{22}^{-1}L_{21} & L_{12} \\ 0 & L_{22} \end{bmatrix}  =  \\[4mm]
&& \Ln^T (x^{*}) \begin{bmatrix} Z_1 & Z_2 \end{bmatrix} 
\begin{bmatrix} L_{11} - L_{12}L_{22}^{-1}L_{21} & L_{12} \\ 0 & L_{22} \end{bmatrix}  \} \\[4mm]
& \subset  \{x^{**} \mid & \Ln^T (x^{**})\hat{Z}\hat{L} = \Ln^T(x^*)\hat{Z}\hat{L} \} = \hat{\mathcal{E}},
\end{array}
\]
since $\hat{Z} = Z_1$ and $\hat{L} = L_{11} - L_{12}L_{22}^{-1}L_{21} $.

For the second statement assume that the full-order network has deficiency zero, i.e., $\ker Z \bigcap \im B =0$, or equivalently $\ker Z \bigcap \im L =0$. Post-multiplication of $L$ with the matrix in (\ref{post}) yields
\[
\ker \begin{bmatrix} Z_1 & Z_2 \end{bmatrix} \bigcap \im \begin{bmatrix} L_{11} - L_{12}L_{22}^{-1}L_{21} & L_{12} \\ 0 & L_{22} \end{bmatrix} =0,
\]
which implies $\ker \hat{Z} \bigcap \im \hat{L} = 0$, i.e., zero-deficiency of the reduced-order network $\hat{\Sigma}$.
\end{proof}

The precise properties of this model reduction approach are currently under investigation. At least the approach seems to have favorable properties for the shortening of linear pathways, leaving out the intermediary complexes. Furthermore, in most cases it seems natural {\it not} to leave out complexes containing boundary chemical species. Note that the model reduction approach results in the reduction of a number of chemical species in the reduced network once these chemical species are only contained in the complexes that are deleted in the reduction procedure.

\subsection{Reversible enzymatic reaction}
As a very simple example, let us apply the above procedure to a standard reversible enzymatic reaction represented as
\[
E + X \rightleftharpoons I \rightleftharpoons E + Y \, ,
\]
where $X, Y$ denote chemical species, $E$ the enzyme, and $I$ the intermediate species (or complex). The three complexes involved are $E + X, I, E+Y$ and we want to leave out the intermediate complex $I$. Since the transposed stoichiometric matrix $S^T$ has full row-rank the Wegscheider conditions (\ref{wegscheider1}) are trivially satisfied, guaranteeing the existence of a thermodynamic equilibrium. Denote for later use the concentrations of the species $X,Y,E$ by $x,y,e$, and their thermodynamic equilibrium concentrations by $x^*,y^*,e^*$. 

Furthermore, denote the balanced reaction constant of the first reversible reaction by $\kappa_X$, and of the second one by $\kappa_Y$. Then the weighted Laplacian matrix of the complex graph is given by
\[
L = \begin{bmatrix}
\kappa_X & - \kappa_X & 0 \\
-\kappa_X & \kappa_X +\kappa_Y & - \kappa_Y \\
0 & -\kappa_Y & \kappa_Y
\end{bmatrix}
\]
Deletion of the intermediate complex $I$ corresponds to taking the Schur complement of $L$ with respect to the second row and column, yielding the reduced weighted Laplacian matrix
\[
\hat{L} = \begin{bmatrix}
\kappa & - \kappa \\
-\kappa & \kappa
\end{bmatrix} \, ,
\]
where the reduced balanced reaction constant $\kappa$ is computed as
\[
\kappa := \frac{\kappa_X \kappa_Y}{\kappa_X + \kappa_Y}
\]
Furthermore, the reduced complex stoichiometric matrix $\hat{Z}$ expressing the remaining complexes $E + X, E + Y$ in terms of the species $X,Y,E$ is given by\footnote{Strictly speaking, the matrix $\hat{Z}$ should be complemented by a zero row corresponding to the chemical species $I$, which is however not present anymore in the remaining complexes.}
\[
\hat{Z} = \begin{bmatrix} 1 & 0 \\ 0 & 1 \\ 1 & 1 \end{bmatrix}
\]
It follows that the reduced dynamics (\ref{reduced}) is given as
\[
\begin{array}{rl}
\begin{bmatrix} \dot{x} \\ \dot{y} \\ \dot{e} \end{bmatrix}  & = -  \begin{bmatrix} 1 & 0 \\ 0 & 1 \\ 1 & 1 \end{bmatrix}
\begin{bmatrix}
\kappa & - \kappa \\
-\kappa & \kappa
\end{bmatrix} \Exp \left( \begin{bmatrix} 1 & 0 & 1 \\ 0 & 1 & 1 \end{bmatrix} 
\begin{bmatrix} \ln \left(\frac{x}{x^*}\right) \\ \ln \left(\frac{y}{y^*}\right) \\ \ln \left(\frac{e}{e^*}\right) \end{bmatrix} \right) \\[4mm]
& =  \begin{bmatrix} -1 \\ 1 \\ 0 \end{bmatrix} \left(\frac{\kappa}{e^* x^*} ex - \frac{\kappa}{e^* y^*}ey \right) ,
\end{array}
\]
which is in mass action kinetics form, and has $(x^*,y^*,e^*)$ as a thermodynamic equilibrium. Note that in the reduced model the total enzyme concentration $e$ is preserved (as to be expected).

It is easily checked that the equations for the $x,y$ and $e$ values of the thermodynamic equilibria of the reduced network are the same as for the original network. On the other hand, in contrast with the thermodynamic equilibrium equations for $i$ of the original model, the reduced model allows any value of $i$ as equilibrium value; in accordance with Proposition \ref{prop:reduction}. Furthermore, both networks are trivially zero-deficient.

\section{Conclusions}\label{sec:Conc}
In this paper we have provided a compact, geometric, formulation of the dynamics of mass action chemical reaction networks possessing a thermodynamic equilibrium. This formulation clearly exhibits both the structure of the complex graph and the stoichiometry, and furthermore admits a direct thermodynamical interpretation. Exploiting this formulation we were able to recover, for this subclass of mass action kinetics chemical reaction networks, some of the results in the fundamental work \cite{HornJackson, Horn, Feinberg, Feinberg1} in a simple and insightful way, without having to rely on the properties of deficiency zero or one. Drawing inspiration from \cite{op, opk, opcircuits}, we have shown how the framework leads to an input-state-output formulation of open chemical reaction networks, and how this may be used for interconnection.
Furthermore, we have indicated how this formulation, in particular its emphasis on the Laplacian matrix of the complex graph, may be used for a systematic model reduction procedure based on the elimination of certain intermediate complexes.

Current research is taking place in a number of directions. We believe the striking similarity with consensus dynamics can be further exploited. Also, the use of our formulation for the analysis of steady states corresponding to non-zero (constant) boundary fluxes is under study, taking into account the possibility of multiple steady states \cite{craciun}; see also the work \cite{Chaves} on input-output stability. The ramifications of the new model reduction procedure are being investigated, including the comparison with classical reduction schemes like (reversible) Michaelis-Menten reaction rates; see e.g. \cite{klipp}. Furthermore, the Laplacian matrix can not only be used for model reduction, but also for other analysis purposes of the complex graph, e.g. for its decompostion \cite{Anderson}. Finally, a challenging question is the extension to {\it biochemical} reaction networks, where the reactions are mostly enzymatic reactions, and the amount of enzymes, through the activity of the gene and protein networks, will depend on e.g. the concentration of a number of chemical species (metabolites). This will lead to the introduction of an additional network structure originating from the regulatory, and possibly competing, feedback loops.

\end{document}